\documentclass[12pt]{article}
\usepackage{graphicx}
\usepackage{amsmath,amsthm,amssymb,enumerate}
\usepackage{euscript,mathrsfs}
\usepackage{color}
\usepackage{dsfont}
\usepackage[left=2cm,right=2cm,top=3.5cm,bottom=3.5cm]{geometry}
\usepackage{color}
\usepackage[framemethod=tikz]{mdframed}
\usepackage{bm}
\allowdisplaybreaks

\usepackage{soul}

\catcode`\@=11 \@addtoreset{equation}{section}

\catcode`\@=12

\allowdisplaybreaks

\newtheorem{Theorem}{Theorem}[section]
\newtheorem{Proposition}[Theorem]{Proposition}
\newtheorem{Lemma}[Theorem]{Lemma}
\newtheorem{Corollary}[Theorem]{Corollary}

\theoremstyle{definition}
\newtheorem{Definition}[Theorem]{Definition}

\newtheorem{Remark}[Theorem]{Remark}

\newcommand{\bTheorem}[1]{
\begin{Theorem} \label{T#1} }
\newcommand{\eT}{\end{Theorem}}

\newcommand{\bProposition}[1]{
\begin{Proposition} \label{P#1}}
\newcommand{\eP}{\end{Proposition}}

\newcommand{\bLemma}[1]{
\begin{Lemma} \label{L#1} }
\newcommand{\eL}{\end{Lemma}}

\newcommand{\bCorollary}[1]{
\begin{Corollary} \label{C#1} }
\newcommand{\eC}{\end{Corollary}}

\newcommand{\bRemark}[1]{
\begin{Remark} \label{R#1} }
\newcommand{\eR}{\end{Remark}}

\newcommand{\bDefinition}[1]{
\begin{Definition} \label{D#1} }
\newcommand{\eD}{\end{Definition}}

\newcommand{\Del}{\Delta_x}

\newcommand{\bfphi}{\boldsymbol{\varphi}}

\newcommand{\bFormula}[1]{
\begin{equation} \label{#1}}
\newcommand{\eF}{\end{equation}}

\newcommand{\Ov}[1]{\overline{#1}}

\newcommand{\aleq}{\stackrel{<}{\sim}}

\newcommand{\ve}{\varepsilon}

\newcommand{\vr}{\varrho}
\newcommand{\vre}{\vr_{\ve}}

\newcommand{\vue}{\vu_{\ep}}

\newcommand{\vu}{\vc{u}}

\newcommand{\vc}[1]{{\bm #1}}

\newcommand{\Div}{{\rm div}_x}
\newcommand{\Grad}{\nabla_x}

\newcommand{\dx}{\,{\rm d} {x}}

\newcommand{\dt}{\,{\rm d} t }

\newcommand{\intO}[1]{\int_{\Omega} #1 \ \dx}

\newcommand{\ep}{\varepsilon}

\newcommand{\br}{ \nonumber \\ }

\newcommand{\mc}{\mathcal}

\def\softd{{\leavevmode\setbox1=\hbox{d}%
          \hbox to 1.05\wd1{d\kern-0.4ex{\char039}\hss}}}
\definecolor{Cgrey}{rgb}{0.85,0.85,0.85}
\definecolor{Cblue}{rgb}{0.50,0.85,0.85}
\definecolor{Cred}{rgb}{1,0,0}
\definecolor{fancy}{rgb}{0.10,0.85,0.10}

\newcommand\Cbox[2]{%
    \newbox\contentbox%
    \newbox\bkgdbox%
    \setbox\contentbox\hbox to \hsize{%
        \vtop{
            \kern\columnsep
            \hbox to \hsize{%
                \kern\columnsep%
                \advance\hsize by -2\columnsep%
                \setlength{\textwidth}{\hsize}%
                \vbox{
                    \parskip=\baselineskip
                    \parindent=0bp
                    #2
                }%
                \kern\columnsep%
            }%
            \kern\columnsep%
        }%
    }%
    \setbox\bkgdbox\vbox{
        \color{#1}
        \hrule width  \wd\contentbox %
               height \ht\contentbox %
               depth  \dp\contentbox
        \color{black}
    }%
    \wd\bkgdbox=0bp%
    \vbox{\hbox to \hsize{\box\bkgdbox\box\contentbox}}%
    \vskip\baselineskip%
}

\mdfdefinestyle{MyFrame}{%
	linecolor=black,
	outerlinewidth=1pt,
	roundcorner=5pt,
	innertopmargin=\baselineskip,
	innerbottommargin=\baselineskip,
	innerrightmargin=10pt,
	innerleftmargin=10pt,
	backgroundcolor=white!20!white}

\date{}

\makeindex
\begin{document}


\title{On the motion of a {nearly incompressible viscous fluid} containing a small rigid body}

\author{Eduard Feireisl
	\thanks{The work of E.F. was partially supported by the
		Czech Sciences Foundation (GA\v CR), Grant Agreement
		21--02411S. The Institute of Mathematics of the Academy of Sciences of
		the Czech Republic is supported by RVO:67985840. A.R and A.Z have been partially supported by the Basque Government through the BERC 2022-2025 program and by the Spanish State Research Agency through BCAM Severo Ochoa excellence accreditation SEV-2017-0718 and through project PID2020-114189RB-I00 funded by Agencia Estatal de Investigación (PID2020-114189RB-I00 / AEI / 10.13039/501100011033). A.Z. was also partially supported  by a grant of the Ministry of Research, Innovation and Digitization, CNCS - UEFISCDI, project number PN-III-P4-PCE-2021-0921, within PNCDI III.} 
		\and Arnab Roy$^1$ \and Arghir Zarnescu$^{1,2,3}$
}

\date{\today}

\maketitle

\bigskip

\centerline{$^*$  Institute of Mathematics of the Academy of Sciences of the Czech Republic}

\centerline{\v Zitn\' a 25, CZ-115 67 Praha 1, Czech Republic}

\centerline{$^1$ BCAM, Basque Center for Applied Mathematics}

\centerline{Mazarredo 14, E48009 Bilbao, Bizkaia, Spain}

\centerline{$^2$IKERBASQUE, Basque Foundation for Science, }

\centerline{Plaza Euskadi 5, 48009 Bilbao, Bizkaia, Spain}

\centerline{$^3$`Simion Stoilow" Institute of the Romanian Academy,}

\centerline{21 Calea Grivi\c{t}ei, 010702 Bucharest, Romania }

\begin{abstract}
	
We consider the motion of a {compressible viscous fluid 
containing a moving rigid body}	confined to a planar domain 
$\Omega \subset R^2$. The main result states that the influence of the body on the fluid is negligible if 
(i) the {diameter} of the body is small and (ii) the fluid is nearly incompressible (the low Mach number regime). 
The specific shape of the body 
as well as the boundary conditions {on the fluid--body interface are irrelevant}
and collisions with the boundary $\partial \Omega$ are allowed. 
The rigid body motion may be enforced externally or governed solely by its interaction with the fluid.

\end{abstract}

\bigskip


\section{Introduction} 
\label{i}

There {is a vast number of recent studies} concerning the motion of a rigid body immersed in/or containing a compressible viscous fluid. {We focus 
on the situation when the body is ``small'' therefore its influence on the fluid motion 
is expected to be negligible. By \emph{small} we mean that the body is contained in a ball 
with a small radius. The problem is mathematically more challenging in the case of 
planar ($2d$) flows, where even small objects may have large capacity.} 

The motion of a small object immersed in an inviscid (Euler) incompressible fluid is studied by Iftimie, Lopes Filho, and Nussenzveig Lopes \cite{MR1974460}. Similar problems \, again in the framework of 
inviscid fluids have been considered by Glass, Lacave, and Sueur \cite{MR3295721}, 
\cite{MR3452278}. The asymptotic behavior of solutions of the incompressible Euler equations in the exterior of a single smooth obstacle when the obstacle becomes very thin tending to curve has been studied by Lacave \cite{MR2542717}.

{
In the context of viscous Newtonian fluids, the flow around a small rigid obstacle 
was studied by Iftimie et al. \cite{MR2244381}.
Lacave \cite{MR2557320} studies the limit of a viscous fluid flow
in the exterior of a thin obstacle shrinking to a curve.} 

{Finally, let us mention results in planar domains, where the body does not influence the flow in the asymptotic limit. Dashti and Robinson \cite{MR2781594} consider 
the viscous fluid-rigid disc system, where the disc is not rotating. Lacave and Takahashi
\cite{MR3595367} consider a single disk moving under the influence of a viscous fluid. They proved convergence
towards the Navier-Stokes equations as the size of the solid tends to zero, its density  is constant and the initial data small. Finally, He and Iftimie \cite{MR3992086} extend the above result to a general shape of the body and to the initial velocities not necessarily small.}
	
{To the best of our knowledge, the problem of negligibility of a small rigid 
	body immersed in a planar viscous \emph{compressible} fluid is completely open. Bravin and  
 Ne\v casov\' a \cite{BraNec} addressed the problem in the $3d$ setting, where 
the capacity of the object in a suitable Sobolev norm is small enough.}

\subsection{Problem formulation} 

 Neglecting completely the possible thermal effects {as well as the external body forces}, we consider the isentropic compressible fluid in the 
low Mach number regime governed by the following system of equations:

\begin{mdframed}[style=MyFrame]
	
	\textsc{Navier--Stokes system.}
	
	\begin{align}
		\partial_t \vr + \Div (\vr \vu) &= 0, \label{i1}\\
		\partial_t (\vr \vu) + \Div (\vr \vu \otimes \vu) + \frac{1}{\ep^{2m}}\Grad p &= \Div \mathbb{S} (\Grad \vu),
		\label{i2} \\
		\mathbb{S}(\Grad \vu) &= \mu \left( \Grad \vu + \Grad^t \vu - \Div \vu \mathbb{I} \right) +
		\lambda \Div \vu \mathbb{I},\ \mu > 0,\ \lambda \geq 0, \label{i3} \\
		p = p(\vr) = a \vr^{\gamma},\ \gamma > 1,\ a>0. \label{i4}
	\end{align}
	
\end{mdframed}

{The fluid is confined to a bounded planar domain $\Omega \subset R^2$ and the momentum equation \eqref{i2} satisfied in} 
\begin{equation} \label{i8}	
	\Omega_{\ep, t} = \Omega \setminus B_{\ep, t},\ t \in (0,T),	
\end{equation} 
{where} 
\begin{equation} \label{i5}
	B_{\ep,t} 	= \left\{ x  \in R^2 \ \Big| \ |x - \vc{h}_{\ve} (t) | \leq \ep \right\},
\end{equation}
\begin{equation} \label{i6}
	\vc{h}_{\ve} \in W^{1,\infty} ([0,T]; R^2),\	\ep |\vc{h}'_\ep (t)| \to 0 \ \mbox{uniformly for a.a.}\ t \in (0,T)\mbox{ as }\ep \to 0.
\end{equation}
{The ball $B_{\ep, t}$ is the part of the plane containing the rigid object 
	at the time $t$. Note carefully that, in general, we do not require $B_{\ep, t} \subset \Omega$. Finally, we impose the no-slip boundary conditions }
\begin{equation} \label{i7}
	\vu|_{\partial \Omega} = 0.
\end{equation}


\subsection{Main results}

Below, we formulate the main hypotheses imposed on the fluid motion. It is convenient to consider the density 
$\vr = \vre$ as well as the velocity $\vu = \vue$ to be defined on the whole physical space $(0,T) \times R^2$. 
Accordingly, we set 
\begin{align} 
	\vr &= \vre (t,x) = \Ov{\vr}  \ - \mbox{a positive constant whenever}\ x \in R^2 \setminus \Omega,\br 
	\vu &= \vue(t,x) = 0 \ \mbox{if} \ x \in R^2 \setminus \Omega.
	\label{i9}
\end{align}

Throughout the whole text, we assume the following:

\begin{itemize}
	
	\item [{\bf (H1)}] 
	\begin{equation}\label{hreg}
	\vc{h}_\ep \in W^{1,\infty} ([0,T]; R^2); 
	\end{equation}
	
	\item[{\bf (H2)}] 
	$(\vre, \vue)$, {$\vre \geq 0$} is a weak renormalized solution of the equation of continuity \eqref{i1}, 
	meaning 
		\begin{align}
		\int_0^T \int_{R^2} \Big[ \vre \partial_t \varphi + \vre \vue \cdot \Grad \varphi \Big] \dx 
		\dt &= - \int_{R^2} \vr_{0, \ep} \varphi (0, \cdot) \ \dx , \br 
		\int_0^T \int_{R^2} \Big[ b(\vre) \partial_t \varphi + b(\vre) \vue \cdot \Grad \varphi  
		+ \left( b(\vre) - b'(\vre) \vre \right) \Div \vue \varphi \Big] & \dx 
	\dt\br &= - \int_{R^2} b(\vr_{\ep,0}) \varphi (0, \cdot) \ \dx  , 
	\label{i10}	
	\end{align}
	for any $\varphi \in C^1_c([0,T) \times R^2)$ and any $b \in C^1[0, \infty)$, $b' \in C_c[0, \infty)$;
	
	\item[{\bf (H3)}] 
$(\vre, \vue)$ is a {weak solution} of the momentum equation \eqref{i2} in the fluid domain 
$\cup_{t \in (0,T)} \Omega_{\ep,t}$, meaning
\begin{equation} \label{i11}
\vue \in L^2(0,T; W^{1,2}_0 (\Omega; R^2)),	
	\end{equation}
and
\begin{align} 
\int_0^T &\intO{ \Big[ \vre \vue \cdot \partial_t \bfphi + \vre \vue \otimes \vue : \Grad \bfphi 
	+ \frac{1}{\ep^{2m}} p(\vre) \Div \bfphi \Big] } \dt \br &= 
\int_0^T \intO{  \mathbb{S} (\Grad \vue) : \Grad \bfphi  } \dt	
- \intO{ \vr_{\ep,0} \vu_{\ep, 0} \cdot \bfphi(0, \cdot) }
	\label{i12}
	\end{align}
for any $\bfphi \in  C^1_c ( \cup_{0 \leq t < T} \Omega_{\ep, t}; R^2) \cap C^1_c([0,T) \times \Omega; R^2)$;

	\item[{\bf (H4)}] The energy inequality 
\begin{align} 
	\int_{\Omega} & \frac{1}{2} \vre |\vue|^2 (\tau, \cdot) \dx + \frac{1}{\ep^{2m}}\int_{\Omega_{\ep, \tau}}   
	\Big( 	P(\vr_{\ep}) - P' (\Ov{\vr} )(\vre - \Ov{\vr} ) - 
	P(\Ov{\vr}) \Big) (\tau, \cdot) \dx \br
	&+ \int_0^\tau \intO{ \mathbb{S}(\Grad \vue) : \Grad \vue } \dt \br &\leq 
	\int_{\Omega}  \frac{1}{2} \vr_{\ep,0} |\vu_{\ep,0}|^2 \dx + \frac{1}{\ep^{2m}}
	\int_{\Omega_{\mathcal{F},\ep,0}}  
	\Big( 	P(\vr_{\ep, 0}) - P' (\Ov{\vr} )(\vr_{\ep,0} - \Ov{\vr} ) - 
	P(\Ov{\vr}) \Big)  \dx
	\label{i13}
\end{align}
holds for a.a. $\tau \in (0,T)$, where $P$ is the pressure potential, 
\[
P(\vr) = \frac{a}{\gamma - 1} \vr^\gamma, \ { \mbox{and}\  \Omega_{\ep, 0} \subset \Omega_{\mathcal{F}, \ep, 0}. }	
\]

\end{itemize}

{In \eqref{i13}, $\Omega_{\mathcal{F},\ep,0 }$ is the fluid domain at the initial time, meaning} 
\[
\Omega_{\mathcal{F}, \ep,0 } \setminus \mathcal{B}_0,\ \mathcal{B}_0 \subset B_{\ep,0} 
\ \mbox{the initial position of the rigid body.}
\]

Our main result reads as follows:

\begin{mdframed}[style=MyFrame]
	
	\begin{Theorem} \label{mT1}
		
	Let $\Omega \subset R^2$ be a bounded domain of class $C^{3}$. Let $(\vre, \vue)_{\ep > 0}$ satisfy the 
	hypotheses {\rm (H1)--(H4)}. In addition, suppose 
	\begin{equation} \label{i14}
	\vr_{\ep,0} \geq 0 \mbox{ a.e. in } \Omega,\quad \frac{1}{\ep^{2m}} {\int_{\Omega_{\mathcal{F},\ep,0}} \Big( P(\vr_{\ep,0}) - P' (\Ov{\vr} )(\vr_{\ep, 0} - \Ov{\vr} ) - 
	P(\Ov{\vr}) \Big) \ \dx} \to 0, 
	\end{equation}
where 
\begin{equation} \label{i14A}
\min \left\{ m; \frac{2 m}{\gamma} \right\} > 3.	
	\end{equation}
	\begin{align}
		\vu_{\ep,0} &\to \vu_0 \ \mbox{weakly in}\ L^2(\Omega; R^2),\ 
	\intO{ \vr_{\ep,0} |\vu_{\ep,0} |^2 } \to \intO{ \Ov{\vr} |\vu_{0} |^2 }\ \mbox{as}\ \ep \to 0 , \br
\ \mbox{where}\ \vu_0 &\in W^{2,\infty}(\Omega), \ \Div \vu_0 = 0,\  \vu_0|_{\partial \Omega} = 0; 	
	\label{i15}	
	\end{align}
	\begin{equation} \label{i16} 
		\ep |\vc{h}'_\ep (t)| \to 0 \ \mbox{uniformly for a.a.}\ t \in (0,T)
		\end{equation}
as $\ep \to 0$.	

Then 
\begin{align} 
\sup_{\tau \in [0,T]} &\| \vre (\tau, \cdot) - \Ov{\vr} {\|_{L^\gamma(\Omega_{\ep, \tau})} } \to 0 \mbox{ with }\gamma \mbox{ as in }\eqref{i4}, \label{i17}\\ 
\vue &\to \vu \ \mbox{in}\ L^2(0,T; W^{1,2}_0 (\Omega; R^2)) 
\label{i18}
\end{align}
as $\ep \to 0$, where $\vu$ is the (unique) classical solution of the incompressible Navier--Stokes 
system	
\begin{align} 
\Div \vu &= 0, \br
\Ov{\vr} \partial_t \vu + \Ov{\vr} \Div (\vu \otimes \vu) + \Grad \Pi &= \mu \Del \vu ,\br 
\vu|_{\partial \Omega} &= 0, \br 
\vu(0, \cdot) &= \vu_0
\label{i19}	
	\end{align}
in $(0,T) \times \Omega$.	

\end{Theorem}

	\end{mdframed}

The hypotheses \eqref{i14}, \eqref{i15} correspond to the \emph{well prepared} data in the low Mach number limit, cf. Masmoudi \cite{MAS}. Moreover, as $\vu_0$ belongs to the class \eqref{i15}, the standard maximal 
regularity theory yields a strong solution of the Navier--Stokes system \eqref{i19}, unique in the class 
\begin{align} 
\vu &\in L^p(0,T; W^{2,p} (\Omega; R^2)),\ 
\partial_t \vu \in L^p(0,T; L^{p} (\Omega; R^2))
, \br 
{\Grad \Pi} &\in L^p(0,T; L^p(\Omega; R^2)),\ 1 \leq p < \infty
\label{i20}
\end{align}
see e.g. Gerhardt \cite{Gerh}, von Wahl \cite{vonWahl}. The solution is classical in $(0,T) \times \Omega$ as a consequence of the interior regularity estimates.

{The hypotheses of Theorem \ref{mT1} are satisfied  if $(\vre , \vue)$ is a weak solution of the fluid--structure interaction problem
	of a single rigid body immersed in a viscous compressible fluid  
	in the sense of \cite{EF64} (see also Desjardins 
	and Esteban \cite{2DEES}) or if the motion of the body is prescribed as in 
	\cite{FKNNS}. A detailed proof is given in Appendix \ref{Ap}.}

The remaining part of the paper is devoted to the proof of Theorem \ref{mT1}. Similarly to the purely incompressible setting studied by He and Iftimie \cite{MR4311107} (cf. also Lacave and Takahashi \cite{MR3595367}), 
the main problem is the rather weak estimate \eqref{i16} that does not allow for a precise identification of 
the limit trajectory of the body. In addition, two new difficulties appear in the compressible regime:

\begin{itemize}
	\item Possible fast oscillations of acoustic (gradient) component of the velocity that cannot be {\it a priori} excluded even for the well prepared data because of the influence of the rigid body.
	
	\item Possible contacts of the body -- intersection of the balls $B_{\ep,t}$--  with the outer boundary 
	$\partial \Omega$.
	
	\end{itemize}

To overcome the above mentioned difficulties, we proceed as follows. In Sections \ref{I}, \ref{II} we identify the system of equations satisfied by the limit velocity $\vu$. Due to the lack of information on 
$\partial_t \vue$, the limit of the convective term as well as the kinetic energy is described in terms 
of the corresponding Young measure. The limit $\vu$ is therefore a generalized dissipative solution of the incompressible Navier--Stokes system in the sense of \cite{AbbFei2}. In particular, we adapt the 
approximation of the test functions introduced by He and Iftimie to the geometry of a bounded domain. Finally, 
in Section \ref{C}, apply the weak--strong uniqueness result proved in \cite{AbbFei2} to conclude that 
the limit is, in fact, a strong solution of the Navier--Stokes system whereas the associated Young--measure reduces to a parametrized family of Dirac masses.

\section{Identifying the limit, the equation of continuity, energy balance}
\label{I}	

It follows from the hypotheses \eqref{i14}, \eqref{i15} that the initial energy on the right--hand side 
of the energy inequality \eqref{i13} is bounded uniformly for $\ep \to 0$. Applying Korn--Poincar\' e 
inequality we get, up to a suitable subsequence,  
\begin{equation} \label{I1}
	\vue \to \vu \ \mbox{weakly in}\ L^2(0, T; W^{1,2}_0 (\Omega; R^2)).
\end{equation}

Next, $\vre$ satisfies the renormalized equation of continuity \eqref{i10}. Moreover, the energy inequality 
\eqref{i13} yields 
\[
\vre \to \Ov{\vr} \ \mbox{in}\ (0,T) \times \Omega \ \mbox{in measure.}
\]
In particular, we may perform the limit in \eqref{i10} obtaining 
\[
b'(\Ov{\vr}) \Ov{\vr} \Div \vu = 0, 
\]
yielding
\begin{equation} \label{I3}
	\Div \vu = 0.
\end{equation}

Finally, using the hypotheses \eqref{i15}, \eqref{i16} and the property of weak lower semi--continuity of convex functionals, we perform the 
the limit in the energy inequality obtaining
	\begin{align} 
	\int_{\Omega} & \frac{1}{2} \Ov{\vr} |\vu|^2 (\tau, \cdot) \dx + \mathfrak{E}(\tau) 
	+ \mu \int_0^\tau \intO{ \Grad \vu : \Grad \vu } \dt  \leq 
	\int_{\Omega}  \frac{1}{2} \Ov{\vr} |\vu_{0}|^2 \ \dx
	\label{I4}
\end{align}
for a.a. $\tau \in (0,T)$. Here, $\mathfrak{C}(\tau) \in L^\infty (0,T)$ is the so called 
\mbox{total energy defect defined as} 
\begin{equation} \label{I4a}
	\mathfrak{E}(\tau) = \liminf_{\ep \to 0} \intO{ \frac{1}{2} \vre |\vue|^2 (\tau, \cdot) } - 
	\intO{ \frac{1}{2} \Ov{\vr} |\vu|^2 (\tau, \cdot) } \geq 0 \ \mbox{for a.a.}\ \tau \in (0,T).
\end{equation}

\section{Identifying the limit, the momentum equation}
\label{II}

The next and more delicate step is to perform the limit $\ep \to 0$ in the momentum equation \eqref{i2}. 
To eliminate the singular pressure term, we consider the test functions 
\begin{equation} \label{I6}
	\bfphi_\ep \in C^1_c ( \cup_{0 \leq t < T} \Omega_{\ep, t}; R^2) \cap 
	C^1_c([0,T) \times \Omega; R^2),\ \Div \bfphi_\ep = 0.	
\end{equation}
Accordingly, the weak formulation \eqref{i12} gives rise to 
\begin{align} 
\int_0^T \intO{ \Big[ \vre \vue \cdot \partial_t \bfphi_\ep + \vre \vue \otimes \vue : \Grad \bfphi_\ep \Big] 
} \dt &= \int_0^T \intO{ \mathbb{S} (\Grad \vue) : \Grad \bfphi_\ep } \dt \br & - 
\intO{\vr_{0,\ep} \vu_{0, \ep} \cdot \bfphi_\ep (0, \cdot) }.
\label{I5}
\end{align}

\subsection{Some useful estimates}
Note that \eqref{I5} is relevant only on the fluid part $\cup_{t \in [0,T]} \Omega_{\ep, t}$, where 
the energy inequality \eqref{i13} yields uniform bounds on the density. This motivates the following decomposition of any measurable functions $v$:
\[
v = [v ]_{\rm ess} + [v ]_{\rm res}, 
\]
where 
\[
[v ]_{\rm ess} = v \mathds{1}_{\frac{1}{2} \Ov{\vr} \leq v \leq 2 \Ov{\vr} }.
\]
Thanks to the energy inequality \eqref{i13}, we get 
\begin{equation} \label{I7}
[ \vre ]_{\rm ess} \vue \ \mbox{bounded in}\ L^\infty(0,T; L^2(\Omega)) \cap 
L^2 (0,T; L^q(\Omega)) \ \mbox{for any}\ 1 \leq q < \infty.	
	\end{equation}
Moreover, by the energy inequality, 
\begin{equation} \label{I8}
	[\vre ]_{\rm ess} \to \Ov{\vr} \ \mbox{in measure in}\ ((0,T) \times \Omega);
\end{equation}
whence we conclude 
\begin{align} 
	[\vre ]_{\rm ess} \vue \to \Ov{\vr} \vu 
	\ &\mbox{weakly -(*) in}\ L^\infty (0,T; L^2(\Omega; R^2)), 
	 \mbox{and weakly in}\ L^2 (0,T;  L^q(\Omega; R^2)) \ \mbox{for any}\ 
	1 \leq q < \infty.
	\label{I9}
	\end{align}
In addition, we also have 
\[
\vre \vue = (\vre - \Ov{\vr}) \vue + \Ov{\vr} \vue, 
\]
where, thanks to the energy inequality \eqref{i13}, 
\begin{align} \label{II9}
	\int_{\Omega_{\ep, \tau} } |\vre - \Ov{\vr} | |\vue| \ \dx \aleq 
	\| \vre(\tau, \cdot) - \Ov{\vr} \|_{(L^{\gamma} + L^2) (\Omega_{\ep, \tau}) } \| \vue \|_{W^{1,2}_0(\Omega; R^2)}
	\aleq \ep^{\min \{ m , \frac{2m}{\gamma} \} } \| \vue(\tau, \cdot) \|_{W^{1,2}_0(\Omega; R^2)}
	\end{align}
for any $\tau \in [0,T]$.

Similarly, 
\begin{align}
[\vre ]_{\rm ess} \vue \otimes \vue \ &\mbox{is bounded in }
L^1 (0,T; L^q(\Omega; R^{d \times d})) \cap L^\infty (0,T; L^1(\Omega; R^{d \times d})) \br
&\mbox{for any}\ 1 \leq q < \infty;
\label{I9a}
\end{align}
whence, by interpolation, 
\begin{equation} \label{I9b}
[\vre ]_{\rm ess} \vue \otimes \vue \to \Ov{\vr \vu \otimes \vu } 
\ \mbox{weakly in}\ L^r ((0,T; L^2(\Omega; R^2)) \ \mbox{for some}\ r > 1.
\end{equation}

The tensor $\Ov{\vr \vu \otimes \vu} \in R^{d \times d}_{\rm sym}$ is positively semi--definite and 
\begin{equation} \label{I9c}
\Ov{\vr \vu \otimes \vu} - \Ov{\vr} \vu \otimes \vu \geq 0. 
\end{equation}
Indeed, for any  $\vc{d}\in R^d$:
\[
\left[ \Ov{\vr \vu \otimes \vu} - \Ov{\vr} \vu \otimes \vu \right]: (\vc{d} \otimes \vc{d}) = 
\lim_{\ep \to 0} | \sqrt{ [ \vre ]_{\rm ess} } \vue \cdot \vc{d} |^2 - 
|\sqrt{\Ov{\vr}} \vu \cdot \vc{d} |^2.
\]	
Thus the desired conclusion \eqref{I9c} follows from \eqref{I1}, \eqref{I8} and weak lower--semicontinuity of convex functions. Finally, as 
\[
[\vre]_{\rm ess} |\vue|^2 \leq \vre |\vue|^2, 
\]
we get 
\begin{equation} \label{I9d}
0 \leq \intO{ {\rm trace} \Big[  \Ov{\vr \vu \otimes \vu} - \Ov{\vr} \vu \otimes \vu \Big] } \leq 
2 \mathfrak{E} ,
\end{equation}	
where $\mathfrak{E}$ is the total energy defect appearing on the left--hand side of the energy inequality \eqref{I4}.

As for the residual components, we deduce from the energy inequality 
\begin{equation} \label{I10}
\int_{\Omega_{\ep,\tau}} [\vre ]^\gamma_{\rm res} (\tau, \cdot) \dx \aleq \ep^{2 m},\ 0 \leq \tau \leq T.
\end{equation}
Consequently, by H\" older's inequality, 
\begin{equation} \label{I11}
	\int_{\Omega_{\ep; \tau}} [\vre ]_{\rm res} |\vue| \dx \aleq 
	\ep^{\frac{2m}{\gamma}} \| \vue (\tau, \cdot) \|_{L^q(\Omega; R^d)},\ \frac{1}{\gamma} + \frac{1}{q} = 1,	
\end{equation}
and, similarly, 
\begin{equation} \label{I12}
	\int_{\Omega_{\ep; \tau}} [\vre ]_{\rm res} |\vue \otimes \vue| \dx \aleq 
	\ep^{\frac{2m}{\gamma}} \| \vue (\tau, \cdot) \|^2_{L^q (\Omega; R^d)} ,\ \frac{1}{\gamma} + \frac{2}{q} = 1	
\end{equation}
for a.a. $\tau \in (0,T)$.
	
\subsection{Constructing a suitable class of test functions}
\label{TF}

Our goal is to approximate a test function 
\[
\bfphi \in C^\infty_c([0, T] \times \Omega; R^2),\ \Div \bfphi = 0, 
\]
by a suitable family of admissible test functions $(\bfphi_\ep)_{\ep > 0}$ in \eqref{I5}.

{The test function are obtained following the construction of 
He and Iftimie \cite{MR3992086,MR4311107}, specifically,}
\[
\widetilde{\bfphi}_\ep = \Grad^\perp (\eta^\ep (x - \vc{h}_\ep (t)) \Psi_\ep ),
\]
with the potential $\Psi_\ep$,
\[
\Grad^{\perp} \Psi_\ep = \bfphi \ \mbox{normalized as}\ \Psi_\ep(t, \vc{h}_\ep (t)) = 0.
\]
The cut-off functions $\eta_\ep$ near the disk $D(\vc{h}_\ep(t),\ep)$ are smooth and satisfy the following properties (see \cite[Lemma 3]{MR3992086}):
\begin{align} \label{I13}
|\eta_\ep |&\leq 1, \eta_\ep(y) = 0 \ \mbox{if}\ |y| \leq \ep, 
\eta_\ep (y) = 1 \ \mbox{if}\ |y| \geq \alpha(\ep) \ep,\\ 
|\Grad \eta_\ep (y) | &\aleq \frac{1}{\ep} \frac{1}{\log (\alpha(\ep))},\ 
|\Grad^2 \eta_\ep (y) | \aleq \frac{1}{\ep^2}.
\label{I13a}
\end{align}
where $\alpha(\ep)$ is chosen in such a way that
\begin{equation} \label{I13b}
\alpha(\ep) \to \infty ,\ \alpha(\ep) \ep (1 + |\vc{h}'_\ep(t) |) \to 0 \ \mbox{as}\ \ep \to 0.
\end{equation}
As shown in \cite[Lemma 5]{MR3992086}, the functions $\widetilde{\bfphi}_\ep$ enjoy the following properties: 
\begin{align}
	\widetilde{\bfphi}_\ep, \ \Grad \widetilde{\bfphi}_\ep  &\in C_c(([0,T] \times R^d) \setminus \cup_{t \in [0,T]} {B}_{\ep, t}),
	\ \partial_t \widetilde{\bfphi}_\ep \in L^\infty((0,T) \times R^2; R^2), \ \label{I16} \\
	{\rm dist}[\vc{h}_\ep(\tau); \partial \Omega] &> \ep \alpha (\ep) \ \Rightarrow\  \widetilde{\bfphi}_\ep (\tau, \cdot)|_{\partial \Omega} = 0, \label{I16a} \\
	\widetilde{\bfphi}_\ep &\to \bfphi \ \mbox{strongly in}\ L^\infty(0,T; W^{1,2}(R^2; R^2)) 
	\ \mbox{as}\ \ep \to 0. \label{I17}
	\end{align}

Unfortunately, the functions $\widetilde{\bfphi}_\ep$ do not vanish on $\partial \Omega$ unless 
${\rm dist}[\vc{h}(t); \partial \Omega] > \ep \alpha (\ep)$. To remedy this, we consider a convex combination 
\[
\bfphi_\ep = \chi_\ep(t) \widetilde{\bfphi}_\ep + (1 - \chi_\ep(t)) \bfphi \ \mbox{for suitable} 
\ 0 \leq \chi_\ep (t) \leq 1,\ \chi_\ep \in W^{1,\infty}(0,T).
\]

First observe that, similarly to $\bfphi_\ep$, 
\[
\| \chi_\ep(t) \widetilde{\bfphi}_\ep + (1 - \chi_\ep) \bfphi \|_{L^\infty(0,T; W^{1,2}(\Omega; R^2))} \aleq 1,
\]
and 
\begin{equation} \label{I20}
\bfphi_\ep - \bfphi = \Big( \chi_\ep(t) \widetilde{\bfphi}_\ep + (1 - \chi_\ep) \bfphi \Big) - \bfphi =
\chi_\ep (\widetilde{\bfphi}_\ep - \bfphi) \to 0 \ \mbox{in}\ L^\infty(0,T; W^{1,2}(\Omega; R^2))
\ \mbox{as}\ \ep \to 0.
\end{equation}

Next, we compute the approximation error in the time derivative 
\[
\partial_t \Big( {\chi}_\ep(t) \widetilde{\bfphi}_\ep + (1 - \chi_\ep) \bfphi \Big) - \partial_t \bfphi 
= \chi_\ep (t) ( \partial_t \widetilde{\bfphi}_\ep - \partial_t \bfphi ) + \chi'_\ep (t) (\widetilde{\bfphi}_\ep  - \bfphi ), 
\]
where the former error term 
\[
 \chi_\ep (t) ( \partial_t \widetilde{\bfphi}_\ep - \partial_t \bfphi ) 
\]
can be controlled in $W^{-1,2}$ exactly as in He and Iftimie \cite{MR4311107} since $\chi$ is independent of $x$. As for the latter,  we have 
\[
\chi_\ep'(t) (\widetilde{\bfphi}_\ep - \bfphi ) = \chi_\ep'(t) \Grad^\perp \Big( [\eta_\ep (x - \vc{h}(t) ) - 1 ]\Psi_\ep \Big)
= \Grad^\perp \left[ \chi_\ep'(t) \Big( [\eta_\ep (x - \vc{h}(t) ) - 1 ]\Psi_\ep \Big) \right],
\]
where, in accordance with \eqref{I13}, 
\begin{equation} \label{I14b}
\| \chi'_\ep(t) [\eta_\ep (x - \vc{h}(t) ) - 1 ]\Psi_\ep \|_{L^2(\Omega)}^2 \aleq 
|\chi'_\ep(t)|^2 \ep^2 \alpha^2 (\ep). 
\end{equation}
Thus if
\begin{equation} \label{I14}
	|\chi'_\ep(t)| \aleq |\vc{h}'_\ep (t) | ,
	\end{equation}
the latter error vanishes in $W^{-1,2}$ for $\ep \to 0$ as a consequence of \eqref{I13b}.

For $\delta > 0$ fixed, let $\bfphi \in C^1([0,T) \times \Omega)$ be given such that 
\begin{equation} \label{I15}
\bfphi (t,x) = 0 \ \mbox{whenever}\ {\rm dist}[x, \partial \Omega] \leq 2 \delta. 
\end{equation} 
Finally, we choose
\[
\chi_\ep (t) = H_{\delta} \Big( {\rm dist}[ \vc{h}_\ep (t); \partial \Omega] \Big) , 
\ 0 \leq H_\delta \leq 1,\ H_\delta (z) = 0 \ \mbox{for}\ z \leq \frac{\delta}{2},\ 
H_\delta(z) = 1 \ \mbox{for}\ z \geq \delta,
\]
where $H_{\delta}$ is a Lipschitz function.
We claim that the test functions 
\[
\bfphi_\ep = \chi_\ep (t) \widetilde{\bfphi}_\ep + (1 - \chi_{\ep}(t) ) \bfphi
\]
vanish both on the boundary $\partial \Omega$ and on the balls $B_{\ep, t}$, $t \in [0,T]$. 
First, by construction, the function 
\[
{\chi}_\ep \widetilde{\bfphi}_\ep 
\]
vanishes on $B_{\ep,t}$ for any $t \in [0,T]$. Moreover, if $\chi_\ep > 0$, then, in view of \eqref{I13b}, 
\[ 
{\rm dist}[ \vc{h}_\ep (t), \partial \Omega] > \frac{\delta}{2} > \ep \alpha(\ep) 
\ \mbox{for}\ \ep \ \mbox{small enough.} 
\]
It follows from \eqref{I16a} that $\chi_\ep \bfphi_\ep|_{\partial \Omega} = 0$.

Second, obviously $(1 - \chi_\ep) \bfphi|_{\partial \Omega} = 0$. Next, if $\chi_{\ep} < 1$, 
we have ${\rm dist}[ \vc{h}_\ep (t); \partial \Omega] < \delta$. Thus, in view of \eqref{I15}, 
$(1 - \chi_\ep) \bfphi (t, \cdot)|_{B_{\ep, t}} = 0$ as soon as $\ep < \delta$.

\subsection{Asymptotic limit}

The function $\bfphi_\ep$ constructed in Section \ref{TF} represents a legitimate test function for 
the momentum balance \eqref{I5}. Our goal is to perform the limit $\ep \to 0$.

\medskip 

\noindent \underline{\bf Step 1: Viscous term.} In view of hypothesis \eqref{i15}, \eqref{I1}, and \eqref{I3}, it follows from 
\eqref{I20} that 
\begin{align}
\int_0^T \intO{ \mathbb{S} (\Grad \vue) : \Grad \bfphi_\ep } \dt  &- 
\intO{\vr_{0,\ep} \vu_{0, \ep} \cdot \bfphi_\ep (0, \cdot) }  \br &\to
\mu \int_0^T \intO{ \Grad \vu : \Grad \bfphi } \dt  - 
\intO{\Ov{\vr} \vu_{0} \cdot \bfphi (0, \cdot) }
\label{A1}
\end{align}
for any $\bfphi \in C^\infty_c ([0,T) \times \Omega; R^d)$, $\Div \bfphi = 0$.

\medskip 

\noindent\underline{\bf Step 2: Convective term.} We can write
\begin{equation*}
\int_0^T \intO{  \vre \vue \otimes \vue : \Grad \bfphi_\ep
} \dt = \int_0^T  \intO{ [\vre]_{\rm ess} \vue \otimes \vue : \Grad \bfphi_\ep } \dt + \int_0^T  \intO{ [\vre]_{\rm res} \vue \otimes \vue : \Grad \bfphi_\ep } \dt
\end{equation*}
 We use \eqref{I9b} to obtain 
\begin{align} 
\int_0^T & \intO{ [\vre]_{\rm ess} \vue \otimes \vue : \Grad \bfphi_\ep } \dt 
\to \int_0^T \intO{ \Ov{\vr} \vu \otimes \vu : \Grad \bfphi } \dt \br
&+ \int_0^T \intO{ \Big( \Ov{\vr \vu \otimes \vu} - \Ov{\vr} \vu \otimes \vu \Big) : \Grad \bfphi } \dt.  
\label{A2}	
	\end{align}

\medskip 
\noindent \underline{\bf Step 3: Time derivative.} Using the same arguments as in \cite{MR4311107} combined with \eqref{I14b}, we get 
\begin{align} 
\intO{ \Ov{\vr} \vue \cdot \partial_t \bfphi_\ep } \aleq \| \vue \|_{W^{1,2}_0(\Omega;R^2)} 
\| \partial_t \bfphi_\ep \|_{W^{-1,2}(\Omega; R^2)}	\to 0 \ \mbox{in}\ L^2(0,T).
\label{A3}	
	\end{align}

\medskip \underline{\bf Step 4: Remaining terms.} The final step is to show
\begin{align} 
\int_0^T \int_{\Omega_{\ep, t}} (\vre - \Ov{\vr}) \vue \cdot \partial_t \bfphi_\ep \ \dx \dt &\to 0, \br
\int_0^T \int_{\Omega_{\ep, t}} [\vre]_{\rm res} \vue \otimes \vue : \Grad \bfphi_\ep \ \dx \dt &\to 0.
\label{A4}	
	\end{align}

A direct manipulation reveals 
\begin{align} 
	\| \Grad \bfphi_\ep \|_{L^\infty((0,T) \times \Omega; R^{2 \times 2})} &\aleq 
	\| \nabla^2 \eta_\ep \|_{L^\infty (R^2)} + 1, \br 
\| \partial_t \bfphi_\ep \|_{L^\infty((0,T) \times \Omega; R^{2 \times 2})} &\aleq (1 + |h_\ep'(t)| )(
\| \nabla^2 \eta_\ep \|_{L^\infty (R^2)} + 1).	
	\label{A5}
	\end{align} 
Consequently, in view of \eqref{I13a} and \eqref{II9}, \eqref{I12}, the desired conclusion \eqref{A4} follows as soon as 
\begin{equation} \label{A6}
\min \left\{ m ; \frac{2m}{\gamma} \right\} > 3.
\end{equation}

\section{Proof of the main result}
\label{C}

Summarizing the results obtained in the preceding section we may infer that limit velocity
\[
\vu \in L^\infty(0,T; L^2(\Omega; R^2)) \cap L^2(0,T; W^{1,2}_0(\Omega; R^2))
\]
solves the following problem: 
\begin{align}
\Div \vu &= 0,\ \vu|_{\partial \Omega} = 0; \br	
\int_0^T \intO{ \Big[ \Ov{\vr} \vu \cdot \partial_t \bfphi + 
\Ov{\vr} \vu \otimes \vu : \Grad \bfphi } \dt &= 
\mu \int_0^T \intO{	\Grad \vu : \Grad \bfphi } \dt  
- \intO{ \Ov{\vr} \vu_0 \cdot \bfphi(0, \cdot) } \br 
&- \int_0^T \intO{ \mathfrak{R} : \Grad \bfphi } \dt	
	\label{C1}
	\end{align}
for any $\bfphi \in C^1_c([0,T) \times \Omega)$;
	\begin{align} 
	\int_{\Omega} & \frac{1}{2} \Ov{\vr} |\vu|^2 (\tau, \cdot) \dx + \mathfrak{E}(\tau) 
	+ \mu \int_0^\tau \intO{ |\Grad \vu | } \dt  \leq 
	\int_{\Omega}  \frac{1}{2} \Ov{\vr} |\vu_{0}|^2 \ \dx
	\label{C2}
\end{align}
for a.a. $\tau \in (0,T)$. Here, the tensor $\mathfrak{R} = \Ov{\vr \vu \otimes \vu} - \Ov{\vr} \vu \otimes \vu 
$ is positively semi--definite and satisfies \eqref{I9d}, specifically 
\begin{equation} \label{C3} 
	0 \leq \intO{ {\rm trace}[ \mathfrak{R} ] } \leq 2 \mathfrak{E} \ \mbox{for a.a.}\ \tau \in (0,T).
\end{equation}
Consequently, the limit function $\vu$ is a dissipative solution of the Navier--Stokes system \eqref{i19} in the sense 
of \cite{AbbFei2}. As the initial velocity is regular, the same problem admits a strong solution 
in the class \eqref{i20}. Thus applying the weak--strong uniqueness result \cite[Theorem 2.6. and Remark 
2.5]{AbbFei2} we conclude that $\vu$ coincides with the strong solution of \eqref{i19}. 

Finally, as the strong solution satisfies the energy equality, its follows from 
\eqref{C2} that $\mathfrak{E}= 0$, and 
\[
\int_0^T \intO{ \mathbb{S}(\Grad \vue) : \Grad \vue } \dt \to 
\mu \int_0^T \intO{ |\Grad \vu|^2 }  
\]
yielding the strong convergence claimed in \eqref{i18}. 

Theorem \ref{mT1} has been proved.

\section{Appendix}

\label{Ap}
	
	Our main result (Theorem \ref{mT1}) is valid whenever $(\vre, \vue)_{\ep > 0}$ satisfy the 
	hypotheses {\rm (H1) -- (H4)} along with the conditions \eqref{i14}--\eqref{i18}.
	These hypotheses (see \eqref{hreg}--\eqref{i13}) are satisfied  if $(\vre , \vue)$ is a weak solution of the fluid--structure interaction problem
	of a single rigid body immersed in a viscous compressible fluid  
	in the sense of \cite{EF64} (see also Desjardins 
	and Esteban \cite{2DEES}) or if the motion of the body is prescribed as in 
	\cite{FKNNS}. Let the rigid body $\mathcal{S}_{\varepsilon}(t)$ be a regular, bounded domain and moving inside $\Omega \subset R^2$. The motion of the rigid body is governed by the balance equations for linear and angular momentum. We assume that the fluid domain $\mathcal{F}_{\varepsilon}(t)=\Omega \setminus \overline{\mathcal{S}_{\varepsilon}(t)}$ is filled with a viscous isentropic compressible fluid. Initially, the domain of the rigid body is given by $\mathcal{S}_{\varepsilon,0}$ included in the ball $B_{\varepsilon,0}$
	and $\mathcal{F}_{\varepsilon,0}$ is the domain of the fluid. Let $h_{\ve}$ be the position of the centre of mass and $\beta_{\ve}$ be the angle of rotation of the rigid body. The solid domain at time $t$  is given by
	\begin{equation*} \mc{S}_{\ve}(t)= h_{\ve}(t)+ \mc{R}_{\beta_{\ve}}(t)\mc{S}_{\ve,0}, \end{equation*} 
	where $\mc{R}_{\beta_{\ve}}$ is the rotation matrix, defined by
	\begin{equation*}
		\mc{R}_{\beta_{\ve}}=\begin{pmatrix}
			\cos \beta_{\ve} & -\sin\beta_{\ve}\\ \sin\beta_{\ve} & \cos\beta_{\ve} 
		\end{pmatrix}.
	\end{equation*}

	The evolution of this fluid-structure system can be described by the following equations 
	\begin{align}
		\frac{\partial {\vre}^{\mc{F}}}{\partial t} + \operatorname{div}({\vre}^{\mc{F}} {\vue}^{\mc{F}}) =0, \quad & \ t\in (0,T),\ x\in \mc{F}_{\varepsilon}(t),\label{mass:comfluid}\\
		\frac{\partial} {\partial t}({\vre}^{\mc{F}}{\vue}^{\mc{F}})+ \operatorname{div}({\rho^{\mc{F}}_{\varepsilon}} {\vue}^{\mc{F}}\otimes {\vue}^{\mc{F}})- \operatorname{div} \mathbb{S}(\nabla_x {\vue}^{\mc{F}})+\frac{1}{\varepsilon^{2m}}\nabla p^{\mc{F}} =0,\quad & \ t\in (0,T),\ x\in \mc{F}_{\varepsilon}(t),\label{momentum:comfluid}\\
		m_{\varepsilon}h_{\varepsilon}''(t)= -\int\limits_{\partial \mc{S}_{\varepsilon}(t)} \left(\mathbb{S}(\nabla_x {\vue}^{\mc{F}}\right) - \frac{1}{\varepsilon^{2m}} p_{\varepsilon}^{\mc{F}}\mathbb{I}) \nu_{\ve}\, d\Gamma,\quad &\mbox{in} \  (0,T), \label{linear momentumcomp:body}\\
		J_{\ve}\beta_{\ve}''(t) = -\int\limits_{\partial \mc{S}_{\ve}(t)}  (\mathbb{S}(\nabla_x {\vue}^{\mc{F}}) - \frac{1}{\varepsilon^{2m}} p_{\ve}^{\mc{F}}\mathbb{I}) \nu_{\ve}\cdot (x-h_{\ve}(t))^{\perp}\, d\Gamma,\quad &\mbox{in} \  (0,T), \label{angular momentumcomp:body}
	\end{align}
	the boundary conditions
	\begin{align}
		{\vue}^{\mc{F}} &= h_{\ve}'(t) +\beta_{\ve}'(t)(x- h_{\ve}(t))^{\perp},  &\mbox{for} \ t \in (0,T),\ x\in \partial \mc{S}_{\ve}(t), \label{boundarycomp-1}\\ 
		{\vue}^{\mc{F}} &= 0,  &\mbox{on} \ (t,x) \in (0,T)\times \partial \Omega, \label{boundarycomp-3}
	\end{align}
	and the initial conditions
	\begin{align}
		{\vre}^{\mc{F}}(0,x)={\vr}_{\mc{F}_{0}}(x),\quad  ({\vre}^{\mc{F}}{\vue}^{\mc{F}})(0,x)=q_{\mc{F}_0}(x), & \quad \forall \ x\in \mc{F}_{\varepsilon,0},\label{initial cond}\\
		h_{\ve}(0)=0,\quad h_{\ve}'(0)=\ell_{0},\quad \beta_{\ve}(0)=0,\quad \beta_{\ve}'(0)=\omega_{0}.\label{initial cond:comp}
	\end{align}
	In the above, the outward unit normal to $\partial \mc{F}_{\ve}(t)$ is denoted by $\nu_{\ve}(t,x)$.
	For all $x=(x_{1}, x_{2})
	\in {R}^{2}$, we denote by $x^{\perp}$, the vector $(
	-x_{2}, x_{1})$. Moreover, the constants $m_{\ve}$ and $J_{\ve}$ are the mass and the moment of inertia of the rigid body.

	We want to state the existence result of the fluid-rigid body interaction system \eqref{mass:comfluid}--\eqref{initial cond:comp}. To do so, we extend the density and  the velocity in the following way:
	\begin{equation}\label{ext:vru}
		\vre (t,x) = \begin{cases}
			{\vre}^{\mc{F}} (t,x),\quad x\in \mc{F}_{\varepsilon}(t), \\
			{\vre}^{\mc{S}}(t,x),\quad x\in \mc{S}_{\varepsilon}(t),\\
			\overline{\vr} ,\quad x\in R^2\setminus\Omega,
		\end{cases} \quad 
		\vue(t,x) = \begin{cases}
			{\vue}^{\mc{F}}(t,x),\quad x\in \mc{F}_{\ve}(t), \\
			h_{\ve}'(t) +\beta_{\ve}'(t)(x- h_{\ve}(t))^{\perp},\quad x \in \mc{S}_{\ve}(t),\\ 0 ,\quad x\in R^2\setminus\Omega.
		\end{cases}
	\end{equation}
	\begin{equation}\label{ext:vru0}
		\vr_{\ep,0}(x) = \begin{cases}
			\vr_{\mc{F}_0} (x),\quad x\in \mc{F}_{\ve,0}, \\
			\vr^{\mc{S}}_{\ve}(0,x),\quad x\in \mc{S}_{\ve,0},\\
			\overline{\vr} ,\quad x\in R^2\setminus\Omega,
		\end{cases} \quad 
		q_{\ep,0}(x) = \begin{cases}
			q_{\mc{F}_0},\quad x\in \mc{F}_{\ve,0}, \\
			\vr^{\mc{S}}_{\ve}(0,x)(\ell_{0} + \omega_{0}\times x),\quad x \in \mc{S}_{\ve,0},\\ 0 ,\quad x\in R^2\setminus\Omega.
		\end{cases}
	\end{equation}
	We have the following existence result for the system \eqref{mass:comfluid}--\eqref{initial cond:comp} following \cite[Theorem 4.1]{EF64}:
	\begin{Theorem}
		Let $\Omega\subset R^2$ be a bounded domain and the pressure $p^{\mc{F}}$ be given by the isentropic constitutive law
		\begin{equation*}
			p^{\mc{F}} = p(\vr^{\mc{F}}) = a (\vr^{\mc{F}})^{\gamma},\ \gamma > 1,\ a>0.
		\end{equation*} Let the initial data $(\vr_0, q_0)$ be defined by \eqref{ext:vru0} satisfying
		\begin{align} \label{init}
			\vr_{{0}} \in L^{\gamma}(\Omega),\quad \vr_{{0}} \geq 0 &\mbox{ a.e. in }\Omega,
			\\ \label{init1}
			q_{\mc{F}_{0}}\mathds{1}_{\{\rho_{\mc{F}_0}=0\}}=0 &\mbox{ a.e. in }\Omega,\quad \dfrac{|q_{\mc{F}_{0}}|^2}{\rho_{\mc{F}_0}}\mathds{1}_{\{\rho_{\mc{F}_0}>0\}}\in L^1(\Omega).
		\end{align}
		Then the system \eqref{mass:comfluid}--\eqref{initial cond:comp} admits a variational solution $(\vr_{\ve},)$ in the following sense:
		\begin{gather}
			{\vre} \geq 0, \quad
			{\vre} \in L^{\infty}(0,T; L^{\gamma}(\Omega)),\quad
			{\vue}\in L^2(0,T; W^{1,2}_{0}(\Omega;R^2)),
			\\
			\vue=h_{\ve}'(t)+\beta_{\ve}'(t) (x-h(t))^{\perp} \mbox{ in }\mc{S}_{\ve}(t),
			\\
			\int\limits_{0}^{T} \int\limits_{{R}^2} \left[{\vre}\frac{\partial \phi}{\partial t}+({\vre}{\vue})\cdot \nabla \phi\right]\, dx\,dt =0, \label{weak density}\\
			\int\limits_{0}^{T} \int\limits_{{R}^2} \left[b({\vre})\frac{\partial \phi}{\partial t}+(b({\vre}){\vue})\cdot \nabla \phi+\left(b({\vre})-b'({\vre}){\vre}\right)\operatorname{div}{\vue}\,\phi\right]\, dx\,dt =0, \label{renormalized}
		\end{gather}
		for any $\phi \in C^1_c([0,T) \times R^2)$ and any $b \in C^1[0, \infty)$, $b' \in C_c[0, \infty)$;
		\begin{equation}\label{continuity weak}
			\int\limits_{0}^{T} \int\limits_{{R}^2} \left[({\vre}{\vue})\cdot \frac{\partial \bfphi}{\partial t}+ ({\vre}{\vue} \otimes {\vue}) : \nabla_x \bfphi +\frac{1}{\varepsilon^{2m}} a{\vre}^{\gamma} \operatorname{div} \bfphi\right]\, dx\, dt 
			=   \int\limits_{0}^{T} \int\limits_{{R}^2}\mathbb{S}(\nabla_x \vue) : \nabla_x\bfphi \, dx\, dt, 
		\end{equation}
		for any $\bfphi \in C_c^{\infty}((0,T)\times \Omega)$, with $\mathbb{D}(\bfphi)=0$ in a neighborhood of $\mc{S}_{\ve}(t)$ where $\mathbb{D}\bfphi=\frac{1}{2}\left( \Grad \bfphi + \Grad^t \bfphi \right)$;

		The following energy inequality holds for a.e. $t \in [0,T]$:
		\begin{align} 
			\int_{\Omega} & \frac{1}{2} \vre |\vue|^2 (\tau, \cdot) \dx + \int_{\Omega}  \frac{1}{\ep^{2m}} 
			\Big( 	P(\vr_{\ep}) - P' (\Ov{\vr} )(\vre - \Ov{\vr} ) - 
			P(\Ov{\vr}) \Big) (\tau, \cdot) \dx + \int_0^\tau \intO{ \mathbb{S}(\Grad \vue) : \Grad \vue } \dt\br
			& \leq 
			\int_{\{\vr_{\ep,0}>0\}}  \frac{1}{2} \frac{ |q_{\ep,0}|^2}{\vr_{\ep,0}} \dx + \frac{1}{\ep^{2m}}
			\int_{\Omega}  
			\Big( 	P(\vr_{\ep, 0}) - P' (\Ov{\vr} )(\vr_{\ep,0} - \Ov{\vr} ) - 
			P(\Ov{\vr}) \Big)  \dx,
			\label{fsi:energy}
		\end{align}
		where $P$ is the pressure potential 
		\[
		P(\vr) = \frac{a}{\gamma - 1} \vr^\gamma. 	
		\]
	\end{Theorem}

	\begin{Remark}
		Let us mention that the specific form of the energy inequality \eqref{i13} follows from 
		\cite[Lemma 3.2]{EF64} and \eqref{fsi:energy}.
	\end{Remark}
	
	We can verify the hypotheses {\rm (H1)--(H4)} and apply Theorem \ref{mT1} under certain conditions to obtain the following result in the framework of fluid-rigid body interaction:
	
	\begin{Theorem} \label{fsi:mT1}
		Let $\Omega \subset R^2$ be a bounded domain of class $C^{3}$ and $(\vr_0, q_0)$  satisfy \eqref{init}--\eqref{init1}. Assume that $\mc{S}_{\ve,0}\subset B_{\ve,0}$,
		\begin{align} \label{fsi:i14}
			&\bullet\quad	\frac{1}{\ep^{2m}} \int_{\Omega_{\ep,0}} \Big( P(\vr_{\ep,0}) - P' (\Ov{\vr} )(\vr_{\ep, 0} - \Ov{\vr} ) - 
			P(\Ov{\vr}) \Big) \ \dx \to 0, \mbox{ where }\min \left\{ m; \frac{2 m}{\gamma} \right\} > 3.\\
			&\bullet\quad		 
			\int_{\{\vr_{\ep,0}>0\}}  \frac{1}{2} \frac{ |q_{\ep,0}|^2}{\vr_{\ep,0}} \dx \to \intO{ \Ov{\vr} |\vu_{0} |^2 }\ \mbox{as}\ \ep \to 0 , \mbox{ where } \vu_0 \in W^{2,\infty}(\Omega), \ \Div \vu_0 = 0,\  \vu_0|_{\partial \Omega} = 0.	
			\label{fsi:i15}	\\
			\label{fsi:i16} 
			&\bullet\quad	\mbox{The mass } m_{\ve}\mbox{ verifies that  }\frac{m_{\ve}}{\ve^2} \to \infty \mbox{ as }\ep \to 0.	
		\end{align}
		Then 
		\begin{align} 
			\sup_{\tau \in [0,T]} &\| \vre (\tau, \cdot) - \Ov{\vr} \|_{(L^2 + L^\gamma)(\Omega)} \to 0, \label{fsi:i17}\\ 
			\vue &\to \vu \ \mbox{in}\ L^2(0,T; W^{1,2}_0 (\Omega; R^2)) \label{fsi:i18}
		\end{align}
		as $\ep \to 0$, where $\vu$ is the (unique) classical solution of the incompressible Navier--Stokes 
		system	
		\begin{align} 
			\Div \vu &= 0, \br
			\Ov{\vr} \partial_t \vu + \Ov{\vr} \Div (\vu \otimes \vu) + \Grad \Pi &= \mu \Del \vu ,\br 
			\vu|_{\partial \Omega} &= 0, \br 
			\vu(0, \cdot) &= \vu_0
			\label{fsi:i19}	
		\end{align}
		in $(0,T) \times \Omega$.	
		
	\end{Theorem}
	
	\begin{Remark}
		We want to point out that as observed by He and Iftimie \cite{MR4311107}, assumption \eqref{i16} holds for the fluid--structure interaction problem if the condition \eqref{fsi:i16} satisfies. Observe that the condition \eqref{fsi:i16} implies $\inf \vr_{\ep}^{\mc{S}} \to \infty$, where $\vr_{\ep}^{\mc{S}}$ is the density of the rigid body immersed in the fluid. 
	\end{Remark}



\end{document}